\documentclass[11pt]{amsart}

\usepackage{amsmath,amssymb,amsthm,mathtools,tikz,graphicx}
\usetikzlibrary{arrows.meta, backgrounds, positioning, calc}

\usepackage[bookmarksopen=true,final]{hyperref}

\usepackage{bm}
\usepackage[nobysame]{amsrefs}

\usepackage{fancyhdr} 
\usepackage[normalem]{ulem}

\usepackage{arydshln}

\renewcommand{\Im}{\mathop{\mathrm{Im}}}

\mathtoolsset{showonlyrefs}

\newcommand{\bbC}{{\mathbb{C}}}

\newcommand{\bbN}{{\mathbb{N}}}

\newcommand{\bbR}{{\mathbb{R}}}

\newcommand{\bbZ}{{\mathbb{Z}}}

\newcommand{\m}{{\mathfrak{m}}}

\newcommand{\calM}{{\mathcal M}}

\newcommand{\bdone}{{\boldsymbol{1}}}

\newcommand{\supp}{\text{\rm{supp}}}

\newcommand{\beq}{\begin{equation}}
	\newcommand{\eeq}{\end{equation}}
\newcommand{\ba}{\begin{align*}}
	\newcommand{\ea}{\end{align*}}

\newcommand{\new}[1]{{\color{blue} 1}}

\allowdisplaybreaks
\numberwithin{equation}{section}

\DeclareFontFamily{U}{mathx}{\hyphenchar\font45}
\DeclareFontShape{U}{mathx}{m}{n}{
	<5> <6> <7> <8> <9> <10>
	<10.95> <12> <14.4> <17.28> <20.74> <24.88>
	mathx10
}{}
\DeclareSymbolFont{mathx}{U}{mathx}{m}{n}
\DeclareMathAccent{\widecheck}{0}{mathx}{"71}

\newtheorem{theorem}{Theorem}
\newtheorem{lemma}[theorem]{Lemma}

\newtheorem{corollary}[theorem]{Corollary}

\newtheorem{proposition}[theorem]{Proposition}

\newtheorem{remark}[theorem]{Remark}
\newtheorem{example}{Example}

\title{Perfectness of partial Christoffel transforms of Nikishin systems}
\date{\today}
\author{Rostyslav Kozhan}
\email{kozhan@math.uu.se}
\address{Uppsala University, Uppsala, Sweden}

\begin{document}
	\maketitle
	\begin{abstract}
		We study multiple orthogonal polynomial systems obtained from a Nikishin system of two measures by applying a Christoffel transform to one of the measures at a point $x=x_0$. We show that if $x_0$ separates the supports of the generating measures, then the modified systems remain perfect.   As an application, we prove that the type I multiple orthogonal polynomials of a Nikishin system have no zeros in the gap between the supports of its generating measures.
	\end{abstract}
	
	\section{Introduction}
	
	Given a positive Borel measure $\mu$ on $\bbR$ with infinite support and finite moments of all orders, there exists a unique sequence of monic orthogonal polynomials $(P_n(x))_{n=0}^\infty$ so that $\deg P_n = n$ and
	\begin{equation}
		\int_\bbR P_n(x) x^j d\mu(x) = 0, \qquad j=0,\ldots,n-1. 
	\end{equation} 
	
	For any point $x_0 \in \bbR$, the Christoffel transform of $\mu$ is the measure denoted by $(x-x_0)\mu$, defined via its action on test functions as
\begin{equation}\label{eq:ChrM10}
	f(x) \mapsto \int_\bbR f(x) (x-x_0) d\mu(x).
\end{equation}
	If $x_0$ lies outside the interior of the convex hull of $\supp \,\mu$, then, up to a sign, $(x-x_0)\mu$ is a positive measure and generates its own unique sequence of monic orthogonal polynomials $(\widehat{P}_n)_{n=0}^\infty$. 
	The mapping $\mu \mapsto (x-x_0)\mu$ and the corresponding polynomial transformation $P_n \mapsto \widehat{P}_n$ are classical operations with numerous applications, most notably in numerical analysis, rational approximation, and discrete integrable systems. 

	In the setting of multiple orthogonal polynomials, orthogonality is distributed across a vector of measures $\bm{\mu} = (\mu_1, \mu_2)$. 	Given a multi-index $(n,m) \in \bbN^2$, a type II multiple orthogonal polynomial is a monic non-zero polynomial $P_{n,m}(x)$ such that $\deg P_{n,m} = n+m$, and
	\begin{align}
		\int P_{n,m}(x) x^k d\mu_1(x) &= 0, \qquad k=0, 1, \dots, n-1, \label{eq:intro_ortho1}\\
		\int P_{n,m}(x) x^k d\mu_2(x) &= 0, \qquad k=0, 1, \dots, m-1. \label{eq:intro_ortho2}
	\end{align}
	The vector of type I multiple orthogonal polynomials $(A_{n,m}, B_{n,m})$ is defined by $\deg A_{n,m}\le n-1$, $\deg B_{n,m}\le m-1$ and
	\begin{equation}\label{eq:moprlIintro}
		\int_{\bbR} A_{n,m}(x) x^k \, d\mu_1(x) + \int_{\bbR} B_{n,m}(x) x^k \, d\mu_2(x) = 
		\begin{cases}
			0, \quad & k = 0,  \dots, n+m-2, \\
			1, & k = n+m-1.
		\end{cases}
	\end{equation}
	(with $A_{0,0} = B_{0,0} \equiv 0$). 
	The system $\bm{\mu}$ is said to be perfect if these polynomials $P_{n,m}(x)$, or equivalently $A_{n,m}(x)$ and $B_{n,m}(x)$, uniquely exist for all multi-indices $(n,m) \in \bbN^2$. See, e.g.,~\cite{Ismail} for more details. 
	
	A prominent class of perfect systems are Nikishin systems~\cite{Nik80}. Given two positive measures $\sigma_1$ and $\sigma_2$ such that the convex hulls $I_1 = \operatorname{conv}(\supp\,\sigma_1)$ and $I_2 = \operatorname{conv}(\supp\,\sigma_2)$ have disjoint interiors, the Nikishin system $\mathcal{N}(\sigma_1, \sigma_2)$ is defined by the vector of measures $\bm{\mu} = (\sigma_1, \mathfrak{m}_{\sigma_2}\sigma_1)$, where 
	\begin{equation}
		\mathfrak{m}_{\sigma_2}(x) = \int_{I_2} \frac{d\sigma_2(t)}{x-t}
	\end{equation}
	is the Markov function of $\sigma_2$. 
	
	Applying the simultaneous Christoffel transform $((x-x_0)\sigma_1, (x-x_0)\mathfrak{m}_{\sigma_2}\sigma_1)$ preserves the Nikishin structure (up to a sign) provided $x_0$ does not lie in the interior of $I_1$. Therefore, the resulting modified system is automatically perfect.

	A much more delicate question is whether applying the Christoffel transform asymmetrically, i.e., modifying only one of the measures $\mu_1$ or $\mu_2$, preserves perfectness. These partial modifications break the Nikishin property. The main result of this paper is that such asymmetrically modified systems 
	are perfect, provided $x_0$ separates the supports of $\sigma_1$ and $\sigma_2$.
	
	In what follows, let $\Delta$ be the closed interval separating $I_1$ and $I_2$ (which reduces to a single point if $I_1$ and $I_2$ touch).
	
	\begin{theorem}\label{thm:1-xm_Intro}
		Let $\mathcal{N}(\sigma_1, \sigma_2) \equiv (\sigma_1, \mathfrak{m}_{\sigma_2} \sigma_1)$ be a Nikishin system. 
		Then the partial Christoffel transformed system $(\sigma_1, (x-x_0) \mathfrak{m}_{\sigma_2} \sigma_1)$ is perfect for any $x_0 \in \Delta$.
	\end{theorem}
	
	\begin{theorem}\label{thm:2-xm_Intro}
		Let $\mathcal{N}(\sigma_1,\sigma_2) \equiv (\sigma_1, \mathfrak{m}_{\sigma_2} \sigma_1)$ be a Nikishin system. 
		Then the partial Christoffel transformed system $((x-x_0)\sigma_1,  \mathfrak{m}_{\sigma_2} \sigma_1)$ is perfect for any $x_0 \in \Delta$.
	\end{theorem}
	
	Neither of the modified systems in Theorems~\ref{thm:1-xm_Intro} and~\ref{thm:2-xm_Intro} is of standard Nikishin type. For Theorem~\ref{thm:1-xm_Intro} (in both cones $n \ge m$ and $n \le m$) and for Theorem~\ref{thm:2-xm_Intro} in the cone $n \le m-2$, we find an algebraic reduction to a standard Nikishin system. For Theorem~\ref{thm:2-xm_Intro} in the cone $n \ge m-2$, such a reduction is not possible; instead, we use the generalized Andr\'eief identity to express the stacked Hankel determinant as a multiple integral and show after symmetrization that the integrand is sign-definite. This multiple-integral approach to normality is inspired by~\cite{CouVA,KNik} and is a common tool in random matrix theory~\cite{Bor,Joh,Kui}.

	
	Beyond generating a new class of perfect systems, these results have direct implications for the study of the zeros of multiple orthogonal polynomials, Toda lattice dynamics, bidiagonal factorizations, and Darboux transformations of Hessenberg matrices.
	
	Indeed, the perfectness of these partially transformed Nikishin systems guarantees the existence of a non-degenerate bidiagonal factorization for the step-line recurrence matrix~\cite{BDFLM,BFM23}, and ensures that discrete Toda-type dynamics driven by sequential asymmetric Christoffel transforms remain well-defined. 

%
	
	As another application, combining Theorems~\ref{thm:1-xm_Intro} and~\ref{thm:2-xm_Intro} with the zero criterion from~\cite{KVInterlacing}, we obtain the following result on the zero-free region for type I polynomials of Nikishin systems.
	
	\begin{corollary}\label{cor:GapZeros}
		Let $\mathcal{N}(\sigma_1, \sigma_2) \equiv (\sigma_1, \mathfrak{m}_{\sigma_2} \sigma_1)$ be a Nikishin system.
		Then $A_{n,m}(x) \ne 0$ on $\Delta$ for all $n \ge 1$, and $B_{n,m}(x) \ne 0$ on $\Delta$ for all $m \ge 1$.
	\end{corollary}
	
	\begin{proof}
		Let $x_0 \in \Delta$. Theorems~\ref{thm:1-xm_Intro} and~\ref{thm:2-xm_Intro} guarantee the perfectness of both partially modified systems at $x_0$.
		As established in~\cite[Sect.~4.2]{KVInterlacing}, the perfectness of $((x-x_0)\sigma_1, \mathfrak{m}_{\sigma_2} \sigma_1)$ is equivalent to the condition that $A_{n,m}(x_0) \ne 0$ for all $(n,m) \in \bbN^2$ with $n \ge 1$. 
		Similarly, the perfectness of $(\sigma_1, (x-x_0)\mathfrak{m}_{\sigma_2} \sigma_1)$ is equivalent to $B_{n,m}(x_0) \ne 0$ for all $(n,m) \in \bbN^2$ with $m \ge 1$. 
	\end{proof}

As numerical simulations indicate, the zeros of $A_{n,m}$ and $B_{n,m}$ can appear on $\bbR$ both to the left and to the right of the gap $\Delta$, as well as in $\bbC\setminus\bbR$. Therefore, if the condition $x_0 \in \Delta$ in Theorems~\ref{thm:1-xm_Intro} and~\ref{thm:2-xm_Intro} is not imposed, one runs the risk of $x_0$ coinciding with one of these ``stray'' zeros. As the preceding proof demonstrates, this would destroy perfectness.
	
		
		\smallskip
	
	Nikishin systems, introduced in 1980 in~\cite{Nik80} (see also~\cite{bookNS}) are a cornerstone of Hermite--Pad\'e approximation theory. Their analytic properties, particularly normality, zero behavior, and asymptotics, have been developed extensively by numerous authors~\cite{Apt99,ALR,BBFL,CouVA,DriSta94,FidIllLop04,FidLop02,FidLop11,FidLop11b,GonLop,KNik,KVInterlacing,Kui,LopMed15}, among many others. Meanwhile, the classical Christoffel transform~\cite{Christoffel}, originally introduced in 1858, has been extended to the multiple orthogonality setting in~\cite{ADMVA,BDFLM,BFM22,BFM23,HueMan,KVChristoffel,KVInterlacing,KVRational,ManRojChr}.

	The paper is organized as follows. In Section~\ref{ss:Prel} we review the necessary preliminaries on Stieltjes transforms and multiple orthogonal polynomials. In Section~\ref{ss:Main} we prove Theorems~\ref{thm:1-xm_Intro} and~\ref{thm:2-xm_Intro} and derive algebraic relations connecting the type~I and type~II polynomials, Hankel determinants, and nearest-neighbor recurrence coefficients of the modified systems to those of standard Nikishin systems. In Section~\ref{ss:Sign}, we determine the signs of the Hankel determinants and nearest-neighbor recurrence coefficients across the multi-index lattice. 
	Finally, in Section~\ref{ss:Examples} we illustrate how classical families such as multiple Laguerre of the first kind and Jacobi--Pi\~{n}eiro polynomials naturally fit within this framework, and present new explicit families of perfect discrete systems constructed from multiple Meixner and Hahn weights with shifted supports.
	
\section{Preliminaries}\label{ss:Prel}
\subsection{Measures and their transforms}

Throughout the paper we use the notation $\bbN = \{k \in \bbZ : k \ge 0\}$, 
$\bbR_+ = \{x \in \bbR : x \ge 0\}$, $\bbR_- = \{x \in \bbR : x \le 0\}$, $\bbC_+ = \{z \in \bbC : \Im z > 0\}$,  $\bbC_- = \{z \in \bbC : \Im z < 0\}$. For any set $S\subseteq\bbR$, 
we denote $\mathring{S}$ to be the interior of $S$ in the topology of $\bbR$, and $\operatorname{conv}(S)$ to be the convex hull of $S$.

Given $I\subseteq\bbR$, let $\mathcal{M}(I)$ be the set of positive Borel measures $\mu$ on $\bbR$ with all moments 
\begin{equation}\label{eq:measure}
	\mu_k = \int_\bbR t^k d\mu(t), \quad k \in \bbN,
\end{equation}
finite and with $\supp\,\mu \subseteq I$ consisting of infinitely many points. We simply write $\calM$ for $\calM(\bbR)$. To any $\mu \in \mathcal{M}$, we associate its Markov function (or Stieltjes transform)
\begin{equation}
	\mathfrak{m}_\mu(z) = \int_{\bbR} \frac{d\mu(t)}{z-t},
\end{equation}
which is analytic in $\bbC \setminus \supp\mu$. 

As in the introduction, for any $\mu \in \mathcal{M}$ and any $z_0 \in \bbC$, the one-step Christoffel transform of $\mu$ is the complex measure denoted by $(t-z_0)\mu$, defined via
\begin{equation}\label{eq:Chr}
	d((t-z_0)\mu)(t) = (t-z_0) d\mu(t).
\end{equation}

More generally, for any weight function $w(t)$ locally integrable with respect to $\mu$, we denote by $w\mu$ the measure defined by $d(w\mu)(t) = w(t) d\mu(t)$.

\begin{lemma}\label{lem:transforms}
	\begin{itemize}
		\item[(i)] For any $z_0 \in \bbC$,
		\begin{equation}\label{eq:m-christoffel}
			\mathfrak{m}_{(t-z_0)\mu}(z) = (z-z_0) \mathfrak{m}_\mu(z) - \mu_0, \quad \text{where } \mu_0 = \int_{\bbR} d\mu.
		\end{equation}
		In particular, $z \mathfrak{m}_\mu(z) = \mathfrak{m}_{t\mu}(z) + \mu_0$.
		
		\item[(ii)] If $z_0 \in \bbR \setminus \supp \, \mu$, then
		\begin{equation}\label{eq:m-geronimus}
			\mathfrak{m}_{(t-z_0)^{-1}\mu}(z)  = 
			\frac{\mathfrak{m}_\mu(z) - \mathfrak{m}_\mu(z_0)}{z - z_0}.
		\end{equation}
		
		\item[(iii)] For any $\mu\in \calM(\bbR)$, there exists the ``stripped measure'' ${\mu}^{(1)}\in\calM(\bbR)$ with $\supp\,{\mu}^{(1)}\subseteq \operatorname{conv}(\supp\,\mu)$ and a constant $b \in \bbR$, 
		such that
		\begin{equation}\label{eq:m-stripping}
			\mathfrak{m}_{\mu^{(1)}}(z) = b+\mu_0^{-1} z - \frac{1}{\mathfrak{m}_\mu(z)}.
		\end{equation}
		
		\item[(iv)] For any $\mu \in\calM(\bbR_-)$, there exists a positive measure ${\mu}^{[1]} \in \calM(\bbR_-)$ 
		such that
		\begin{equation}\label{eq:m-tilde}
			\mathfrak{m}_{{\mu}^{[1]}}(z) = \frac{1}{z\mathfrak{m}_\mu(z)} - \mu_0^{-1}.
		\end{equation}
	\end{itemize}
\end{lemma}
\begin{proof}
	\textit{(i)}, \textit{(ii)}, and \textit{(iii)} are standard results in the theory of orthogonal polynomials, see, e.g., \cite{Ismail,SimonL2}.
	
	
	
	For \textit{(iv)}, we have $\Im\big(z\mathfrak{m}_\mu(z)\big) = \int_{\bbR_-} \frac{-t \Im z}{|z-t|^2} d\mu(t) > 0$ for $z \in \bbC_+$. Thus, $z\mathfrak{m}_\mu(z)$ is a Herglotz function, which implies that $G(z) = - (z\mathfrak{m}_\mu(z))^{-1} + c$ is also Herglotz. By the Herglotz representation theorem, $G$ must be $-\mathfrak{m}_{{\mu}^{[1]}}$ for some positive Borel measure ${\mu}^{[1]}$ up to a linear term. The linear term vanishes since $G(z) = O(z^{-1})$ as $z \to \infty$, and expanding $G(z)$ at $\infty$ shows that $\mu^{[1]}$ has finite moments of all orders. Finally, $\supp\,{\mu}^{[1]}\subseteq \bbR_-$ since $z\mathfrak{m}_\mu(z) > 0$ on $(0,\infty)$, so $G(z)$ is analytic and real-valued on $(0,\infty)$.
	
	An alternative way of constructing ${\mu}^{[1]}$ is to define ${\mu}^{[1]} = -t^{-1}\mu^{(1)} + k \,\delta_0$, where $\mu^{(1)}$ is the measure from \textit{(iii)} and $k = b - \mathfrak{m}_{\mu^{(1)}}(0) = (\int_{\bbR_-} |t|^{-1} d\mu(t))^{-1} \ge 0$ (with $k=0$ if the integral diverges). Note that $\int_{\bbR_-} (-t)^{-1} d\mu^{(1)}(t) = b - k < \infty$ by taking $z \to 0^+$ in~\eqref{eq:m-stripping}, so ${\mu}^{[1]} \in \calM(\bbR_-)$ is well-defined. One says that ${\mu}^{[1]}$ is a generalized Geronimus transform of $\mu^{(1)}$, that is, an inverse of the Christoffel transform, since $\mu^{(1)} = - t \mu^{[1]}$ holds.
\end{proof}

\subsection{Multiple orthogonal polynomials}

We provide a brief review of the theory of multiple orthogonal polynomials, for more details, see~\cite{Ismail}. We restrict ourselves to the case of two measures. To this end, consider a system $\bm{\mu} = (\mu_1, \mu_2) \in \mathcal{M}^2$.  

Recall the definition of the type I and type II multiple orthogonal polynomials as in the introduction.

%

The existence and uniqueness of either type I or type II polynomials is equivalent to the condition $H_{n,m}[\bm{\mu}] \ne 0$, where $H_{n,m}[\bm{\mu}]$ is the stacked Hankel $(n+m)\times (n+m)$ determinant constructed from the moments $\mu_{1,k} = \int_{\bbR} t^k d\mu_1(t)$ and $\mu_{2,k} = \int_{\bbR} t^k d\mu_2(t)$:
\begin{equation}\label{eq:Hankel}
	H_{n,m}[\bm{\mu}] = 
	\det \left(\begin{array}{c c c c}
		\mu_{1,0} & \mu_{1,1} & \cdots & \mu_{1,n+m-1} \\
		\mu_{1,1} & \mu_{1,2} & \cdots & \mu_{1,n+m} \\
		\vdots & \vdots & \ddots & \vdots \\
		\mu_{1,n-1} & \mu_{1,n} & \cdots & \mu_{1,2n+m-2} 
		\\[0.4em]
		\hdashline[0.5pt/2pt]
		\\[-1.2em]
		\mu_{2,0} & \mu_{2,1} & \cdots & \mu_{2,n+m-1} \\
		\mu_{2,1} & \mu_{2,2} & \cdots &  \mu_{2,n+m} \\
		\vdots & \vdots & \ddots & \vdots \\
		\mu_{2,m-1} & \mu_{2,m} & \cdots & \mu_{2,n+2m-2} \\
	\end{array}\right),
\end{equation}
with the convention $H_{0,0}[\bm{\mu}] = 1$.

If $H_{n,m}[\bm{\mu}]$ is non-zero, we say that the index $(n,m)$ is normal with respect to $\bm{\mu}$. The system $\bm{\mu}$ is called perfect if every $(n,m) \in \bbN^2$ is normal.



\subsection{Nearest neighbor recurrence relations}

When the system $\bm{\mu}$ is perfect, the monic type II multiple orthogonal polynomials satisfy a system of nearest neighbor recurrence relations~\cite{NNRR}
\begin{align}
	x P_{n,m}(x) &= P_{n+1,m}(x) + b_{n,m;1} P_{n,m}(x) + a_{n,m;1} P_{n-1,m}(x) + a_{n,m;2} P_{n,m-1}(x), \label{eq:nnrr-1} \\
	x P_{n,m}(x) &= P_{n,m+1}(x) + b_{n,m;2} P_{n,m}(x) + a_{n,m;1} P_{n-1,m}(x) + a_{n,m;2} P_{n,m-1}(x), \label{eq:nnrr-2}
\end{align}
for all $(n,m)\in\bbN^2$ with the initial conditions $P_{0,0}(x) \equiv 1$ and $a_{0,m;1}=0$, $a_{n,0;2}=0$.  

The recurrence coefficients are linked to the stacked Hankel determinants via
\begin{equation}\label{eq:a-coeffs}
	a_{n,m;1} = \frac{H_{n+1,m}[\bm{\mu}] H_{n-1,m}[\bm{\mu}]}{\left(H_{n,m}[\bm{\mu}]\right)^2}, \qquad 
	a_{n,m;2} = \frac{H_{n,m+1}[\bm{\mu}] H_{n,m-1}[\bm{\mu}]}{\left(H_{n,m}[\bm{\mu}]\right)^2}.
\end{equation}
and
\begin{equation}\label{eq:b-diff}
	b_{n,m;2} - b_{n,m;1} = 
	\frac{H_{n+1,m+1}[\bm{\mu}] H_{n,m}[\bm{\mu}]}{H_{n+1,m}[\bm{\mu}] H_{n,m+1}[\bm{\mu}]}.
\end{equation}
The sign of the coefficients~\eqref{eq:b-diff} are often of interest as they govern the direction of zero interlacings in Angelesco, AT, and Nikishin systems.

\subsection{Chebyshev and AT properties}

Let $I\subseteq \bbR$ be an open or closed interval. 
A collection $(u_j(t))_{j=1}^N$ of continuous real-valued functions on 
$I$ is called a Chebyshev system on $I$ if the determinant
\begin{equation}\label{eq:Chebyshev}
	U_{\bm{x}}(u_1, \ldots, u_N)=\det
	\begin{pmatrix}
		u_1(x_1) & u_1(x_2) & \cdots & u_1(x_N) \\
		u_2(x_1) & u_2(x_2) & \cdots & u_2(x_N) \\
		\vdots & \vdots & \ddots & \vdots \\
		u_N(x_1) & u_N(x_2) & \cdots & u_N(x_N) \\
	\end{pmatrix} 
\end{equation}
is non-zero for any $x_1<x_2<\ldots <x_N$ on $I$.

In particular, given $n,m\in\bbN$, $N=n+m$, and $\bm{t}\in\bbR^m$, $\bm{x}\in\bbR^N$ with $x_k \ne t_j$, we denote the Cauchy--Vandermonde determinant by
\begin{multline}    \label{eq:C}
	C(\bm{x};\bm{t}) = U_{\bm{x}} \Big( 1, x, \ldots, x^{n-1} , \frac{1}{x-t_1}, \ldots, \frac{1}{x-t_{m}}\Big)
	\\
	=
	(-1)^{S(n,m)} \prod_{1\le j<k\le N} (x_k - x_j) \prod_{1\le j<k\le m} (t_k-t_j)  \prod_{j=1}^{m}\prod_{k=1}^{N}(x_k-t_j)^{-1},
\end{multline}
where
\begin{equation}
	S(n,m) = nm + \frac{m(m-1)}{2}. 
\end{equation}
This identity is well-known: one can subtract the first column from each of the other columns, extract a factor of $(x_k-x_1)$ from column $k$ ($k=2,\dots,N$) and a factor of $-(x_1-t_j)^{-1}$ from the last $m$ rows, and then proceed by induction. See~\cite{GMM} for more general determinants of this type.

Let $d\mu_1(x)=w_1(x)d\mu(x)$ and $d\mu_2(x)=w_2(x)d\mu(x)$ be measures with continuous real-valued weights $w_1, w_2$ with respect to some $\mu\in\mathcal{M}(I)$. Applying the Andr\'eief identity to~\eqref{eq:Hankel}, we obtain the integral representation~\cite{Joh,Kui}:
\begin{multline}\label{eq:HankelIntegralAT}
	H_{n,m}[\bm\mu]
	=
	\tfrac{1}{N!} \int_{\bm{x}\in I^N}    
	U_{\bm{x}} \Big(w_1, xw_1,\ldots,x^{n-1}w_1,w_2,xw_2,\ldots,x^{m-1}w_2 \Big)
	\\
	\times \prod_{1 \le j<k \le N}(x_k-x_j) \, \prod_{j=1}^N d\mu(x_j),
\end{multline}
where $N=n+m$.

We say that $(\mu_1,\mu_2)$ is an AT system on $I$ for the index $\bm{n} = (n,m)\in\bbN^2$ if the functions
\begin{equation}\label{eq:AT}
	\{w_1, xw_1,\ldots,x^{n-1}w_1,w_2,xw_2,\ldots,x^{m-1}w_2\}
\end{equation}
form a Chebyshev system on the open interval $\mathring{I}$. In that case, the integrand in~\eqref{eq:HankelIntegralAT} maintains a strict sign on the simplex $x_1 < \dots < x_N$, which immediately implies $H_{n,m}[\bm{\mu}] \ne 0$.

\subsection{Nikishin systems}



Let $\sigma_1 \in \calM(I_1), \sigma_2 \in \mathcal{M}(I_2)$ where
$I_1$ and $I_2$ are two closed intervals such that
\begin{equation}
	\mathring{I}_1 \cap \mathring{I}_{2} = \varnothing.
\end{equation}
We define the corresponding Nikishin system of two measures as
\begin{equation}
	\mathcal{N}(\sigma_1, \sigma_2) = (\sigma_1, \m_{\sigma_2} \sigma_1).
\end{equation}
A fundamental property of Nikishin systems is that they are AT systems on $I_1$ for any multi-index $(n,m) \in \bbN^2$~\cite{Nik80,DriSta94,BBFL,FidIllLop04,FidLop02,FidLop11,FidLop11b,KNik}. Consequently, every Nikishin system is perfect.


If $I_1$ and $I_2$ share a common endpoint, we implicitly assume that the behavior of $\sigma_1$ and $\sigma_2$ at this point is such that the resulting measure $\m_{\sigma_2} \sigma_1$ remains finite. 

For Nikishin systems,~\eqref{eq:HankelIntegralAT} becomes
\begin{multline}\label{eq:HankelIntegralNik}
	H_{n,m}\big[\mathcal{N}(\sigma_1, \sigma_2)\big]
	\\=
	\tfrac{1}{N!} \int_{\bm{x}\in I_1^N}   
	U_{\bm{x}} \Big(1,x,\ldots,x^{n-1},
	\mathfrak{m}_{\sigma_2}(x),x \mathfrak{m}_{\sigma_2}(x),\ldots,x^{m-1} \mathfrak{m}_{\sigma_2}(x)\Big)
	\\
	\times \prod_{1 \le j<k \le N}(x_k-x_j) \, \prod_{j=1}^N d\sigma_1(x_j).
\end{multline}

The following result is well-known to the experts. We provide the proof for completeness.

\begin{proposition}\label{prop:nnrr-signs}
	Let $\bm{\mu} = \mathcal{N}(\sigma_1, \sigma_2)\equiv (\sigma_1, \m_{\sigma_2}\sigma_1)$ be a Nikishin system. 
	
	Define 
	\begin{equation}\label{eq:flippingNikishin}
		\widecheck{\sigma}_1 = \mathfrak{m}_{\sigma_2}\sigma_1, \qquad  \widecheck{\sigma}_2 = \sigma_2^{(1)},
	\end{equation}
	 where $\sigma_2^{(1)}$ is the ``stripped measure'' of $\sigma_2$ 
	as defined in Lemma~\ref{lem:transforms}\textit{(iii)}. Define the ``reversed'' Nikishin system $\widecheck{\bm{\mu}} = \mathcal{N}(\widecheck{\sigma}_1, \widecheck{\sigma}_2) \equiv (\widecheck{\sigma}_1,\mathfrak{m}_{\widecheck{\sigma}_2} \widecheck{\sigma}_1)$.
	
	Denote $P_{n,m}$ and $\widecheck{P}_{n,m}$ the type II multiple orthogonal polynomials for $\bm{\mu}$ and $\widecheck{\bm{\mu}}$, respectively, and similarly for the type I polynomials $A_{n,m}$ and $B_{n,m}$ and Hankel determinants $H_{n,m}$.
	
	\begin{itemize}
		\item[(i)] For any $n \le m-1$:
		\begin{equation}
			H_{n,m} = (-1)^{n(m+1)} \widecheck{H}_{m,n}.
		\end{equation}
		
		\item[(ii)] For any $n \le m-1$:
		\begin{equation}
			P_{n,m}(x) = \widecheck{P}_{m,n}(x).
		\end{equation}
		
		\item[(iii)] For $n \le m-1$:
		\begin{equation}
			\begin{aligned}
				A_{n,m}(x) &= -\widecheck{B}_{m,n}(x), \\
				B_{n,m}(x) &= \widecheck{A}_{m,n}(x) + (cx+b)\widecheck{B}_{m,n}(x),
			\end{aligned}
		\end{equation}
		where $c = (\int d\sigma_2)^{-1}$ and $b$ is from Lemma~\ref{lem:transforms}\textit{(iii)} for $\sigma_2$.
	\end{itemize}
\end{proposition}

\begin{proof}
	To prove \textit{(i)}, we apply the Andreief identity~\eqref{eq:HankelIntegralNik}.  Factoring $\mathfrak{m}_{\sigma_2}(x_k)$ out of each column of $U_{\bm{x}}(1, \dots, x^{n-1}, \mathfrak{m}_{\sigma_2}, \dots, x^{m-1}\mathfrak{m}_{\sigma_2})$ leaves us with $U_{\bm{x}}(1/\mathfrak{m}_{\sigma_2}, \dots, x^{n-1}/\mathfrak{m}_{\sigma_2}, 1, \dots, x^{m-1})$ integrated against $d\widecheck{\sigma}_1(x_k) = \mathfrak{m}_{\sigma_2}(x_k) d\sigma_1(x_k)$. Using $1/\mathfrak{m}_{\sigma_2}(x) = cx+b - \mathfrak{m}_{\widecheck{\sigma}_2}(x)$ and $n\le m-1$, we can reduce $U_{\bm{x}}$ via elementary row operations to $U_{\bm{x}}(-\mathfrak{m}_{\widecheck{\sigma}_2}, \dots, -x^{n-1}\mathfrak{m}_{\widecheck{\sigma}_2}, 1, \dots, x^{m-1})$.
	Factoring out the minus sign from the top $n$ rows and swapping the top block of $n$ rows with the bottom block of $m$ rows introduces a sign $(-1)^{n+nm}$, with the resulting integral being $\widecheck{H}_{m,n}$ for the system $\widecheck{\bm{\mu}}$. This proves \textit{(i)}.

	 By definition~\eqref{eq:intro_ortho1}--\eqref{eq:intro_ortho2}, $P_{n,m}(x)$ of degree $n+m$ satisfies 
	 \begin{alignat}{3}
	 	\label{eq:PII-1}
	 	& \int_{\bbR} P_{n,m}(x) x^k d\sigma_1(x) = 0, \qquad && k=0, \dots, n-1, 
	 	\\
	 	\label{eq:PII-2}
	 	& \int_{\bbR} P_{n,m}(x) x^j \mathfrak{m}_{\sigma_2}(x) d\sigma_1(x) = 0, \qquad &&j=0, \dots, m-1. 
	 \end{alignat}
	 To show $P_{n,m} = \widecheck{P}_{m,n}$, we must verify that $P_{n,m}$ satisfies the orthogonality conditions for $\widecheck{\bm{\mu}}$ at index $(m,n)$. 
	 The first set of conditions for $\widecheck{P}_{m,n}$ requires orthogonality against $x^j d\widecheck{\sigma}_1(x) = x^j \mathfrak{m}_{\sigma_2}(x) d\sigma_1(x)$ for $j=0, \dots, m-1$. These coincide with \eqref{eq:PII-2}. 
	 
	 The second set of conditions requires orthogonality against $x^k \mathfrak{m}_{\widecheck{\sigma}_2}(x) d\widecheck{\sigma}_1(x)$ for $k=0, \dots, n-1$. Using the definition of $\widecheck{\sigma}_j$'s and~\eqref{eq:m-stripping}, the integrand becomes:
	 \begin{multline}
	 	P_{n,m}(x) x^k \big(cx+b - (\mathfrak{m}_{\sigma_2}(x))^{-1}\big) \mathfrak{m}_{\sigma_2}(x) d\sigma_1(x)
	 	\\ = P_{n,m}(x) \big[ x^k(cx+b)\mathfrak{m}_{\sigma_2}(x) - x^k \big] d\sigma_1(x).
	 \end{multline}
	 Integrating this expression yields two terms. 
	 The first term is zero by~\eqref{eq:PII-2} and $k+1\le n \le m-1$. The second term is zero by~\eqref{eq:PII-1} and $k\le n-1$. Since both $P_{n,m}$ and $\widecheck{P}_{m,n}$ are monic polynomials of degree $n+m$ satisfying the same orthogonality conditions, normality implies  $P_{n,m} = \widecheck{P}_{m,n}$.
	
	To prove \textit{(iii)}, the type I vector $(A_{n,m}, B_{n,m})$ for $\bm{\mu}$ is characterized by the degree bounds $\deg A_{n,m} \le n-1$, $\deg B_{n,m} \le m-1$, and the orthogonality conditions:
	\begin{equation}
		\int_{\bbR} \big(A_{n,m}(x) + B_{n,m}(x) \mathfrak{m}_{\sigma_2}(x)\big) x^k d\sigma_1(x) = \delta_{k, n+m-1}, \quad k=0, \dots, n+m-1.
	\end{equation}
	Using the definition of $\widecheck{\sigma}_j$'s and~\eqref{eq:m-stripping}, we obtain
\begin{multline}
	\int_\bbR \Big( \big[B_{n,m}(x) + (cx+b)A_{n,m}(x)\big] - A_{n,m}(x) \mathfrak{m}_{\widecheck{\sigma}_2}(x) \Big) x^k d\widecheck{\sigma}_1(x)
	= \delta_{k, n+m-1}, 
	\\
	k=0, \dots, n+m-1.
\end{multline}
	Since $n\le m-1$, we have $\deg \big(B_{n,m}(x) + (cx+b)A_{n,m}(x)\big) \le m-1$ and $\deg (-A_{n,m}(x)) \le n-1$. 
	These polynomials match the type I vector $(\widecheck{A}, \widecheck{B})$ for $\widecheck{\bm{\mu}}$ at the reversed index $(m,n)$. By normality we obtain \textit{(iii)}. 
%
\end{proof}


\section{Partial Christoffel Transformed Nikishin Systems}\label{ss:Main}

In this section, we investigate systems obtained from a classical Nikishin system by applying the one-step Christoffel transform~\eqref{eq:Chr} to one of the measures.
We establish that these transformed systems are perfect, provided the shift point $x_0$ separates the convex hulls $I_1$ and $I_2$ of the supports of the generating measures $\sigma_1$ and $\sigma_2$.

Before proceeding, a brief methodological remark is in order. While general perfectness criteria for Christoffel, Geronimus, and Uvarov transforms are known~\cite{KVRational}, applying them here requires \textit{a priori} knowledge of the type I zero locations for arbitrary Nikishin systems. Since these are unknown, we instead prove perfectness directly by constructing algebraic mappings to standard Nikishin systems and applying the Andréief identity. Reversing the logic, this constructive proof then yields the zero-free regions as a natural consequence (Corollary~\ref{cor:GapZeros}).


By translation and reflection, we can reduce the general setting to the canonical case where $I_1 \subseteq \bbR_+$, $I_2 \subseteq \bbR_-$, and $x_0=0$.

\subsection{System $(\sigma_1, x \mathfrak{m}_{\sigma_2}\sigma_1)$}

Here, we analyze the modified system $\bm{\mu} = (\sigma_1, x \mathfrak{m}_{\sigma_2} \sigma_1)$ with $\supp\,\sigma_1 \subseteq \bbR_+$ and $\supp\,\sigma_2 \subseteq \bbR_-$. We will show that the properties of this system are inherited from two standard Nikishin systems, which are coupled across the diagonal multi-indices $(n,m)$ with $n=m$.

We denote the first Nikishin system by $\widehat{\bm{\mu}}= \mathcal{N}(\widehat{\sigma}_1, \widehat{\sigma}_2) \equiv (\widehat{\sigma}_1, \mathfrak{m}_{\widehat{\sigma}_2} \widehat{\sigma}_1)$, where 
\begin{equation}\label{eq:hatsigma_def}
	\widehat{\sigma}_1 = \sigma_1, \qquad \widehat{\sigma}_2 = -t \sigma_2.
\end{equation}

The second Nikishin system is defined by $\widetilde{\bm{\mu}}=\mathcal{N}(\widetilde{\sigma}_1, \widetilde{\sigma}_2) \equiv (\widetilde{\sigma}_1, \mathfrak{m}_{\widetilde{\sigma}_2} \widetilde{\sigma}_1)$, where
\begin{equation}\label{eq:tildesigma_def}
	\widetilde{\sigma}_1 = x\mathfrak{m}_{\sigma_2}\sigma_1, \qquad \widetilde{\sigma}_2 = \sigma_2^{[1]}
\end{equation}
(it is illuminating to compare this with~\eqref{eq:flippingNikishin}), where $\sigma_2^{[1]}$ is the measure from Lemma~\ref{lem:transforms}\textit{(iv)}. From this we have
\begin{equation}\label{eq:sigma2tilde}
	\mathfrak{m}_{\widetilde{\sigma}_2}(x) = \frac{1}{x\mathfrak{m}_{\sigma_2}(x)} - \frac{1}{c},
\end{equation}
where $c = \int_{\bbR_-} d\sigma_2$.

Note that all four measures $\widehat{\sigma}_1, \widehat{\sigma}_2, \widetilde{\sigma}_1, \widetilde{\sigma}_2$ are positive with $\widehat{\sigma}_1$ and $\widetilde{\sigma}_1$ supported on $\bbR_+$ and $\widehat{\sigma}_2, \widetilde{\sigma}_2$ supported on $\bbR_-$. 
Crucially, if the supports share the origin ($0 \in I_1 \cap I_2$), both systems automatically satisfy the Nikishin integrability conditions at the boundary. After substitution, we get
\begin{align}
	d(\mathfrak{m}_{\widehat{\sigma}_2}\widehat{\sigma}_1)(x) &= \big(c - x\mathfrak{m}_{\sigma_2}(x)\big) d\sigma_1(x), \\
	d(\mathfrak{m}_{\widetilde{\sigma}_2}\widetilde{\sigma}_1)(x) &= \big(1 - c^{-1} x\mathfrak{m}_{\sigma_2}(x)\big) d\sigma_1(x),
\end{align}
which are easily seen to be integrable.
 Thus, both $\widehat{\bm\mu}$ and $\widetilde{\bm\mu}$ are indeed well-defined standard Nikishin systems. As such, both $\widehat{\bm\mu}$ and $\widetilde{\bm\mu}$ are perfect, which we use in the proof below. 
 
Let $H_{n,m}, P_{n,m}, A_{n,m}, B_{n,m}$ denote the stacked Hankel determinants and multiple orthogonal polynomials for the modified system $\bm{\mu}$, with analogous notations $\widehat{H}_{n,m}, \widehat{P}_{n,m}, \widehat{A}_{n,m}, \widehat{B}_{n,m}$ and $\widetilde{H}_{n,m}, \widetilde{P}_{n,m}, \widetilde{A}_{n,m}, \widetilde{B}_{n,m}$  for the systems $\widehat{\bm{\mu}}$ and $\widetilde{\bm{\mu}}$, respectively.
 We can now prove the following result.
 
\begin{theorem}\label{thm:XM1_Properties}
	Let $\bm{\mu} = (\sigma_1, x \mathfrak{m}_{\sigma_2} \sigma_1)$ with $\supp\,\sigma_1 \subseteq \bbR_+$ and $\supp\,\sigma_2 \subseteq \bbR_-$, and define $\widehat{\bm{\mu}}$ and $\widetilde{\bm{\mu}}$ to be the Nikishin systems defined in~\eqref{eq:hatsigma_def} and~\eqref{eq:tildesigma_def}.
	\begin{itemize}
		\item[(i)] For any $(n,m)\in\bbN^2$, the Hankel determinants are related via:
		\begin{equation}
			H_{n,m}
			= 
			\begin{cases}
				(-1)^{m} \widehat{H}_{n,m}, & \mbox{if } n \ge m, \\
				(-1)^{nm} \widetilde{H}_{m,n}, & \mbox{if } n \le m.
			\end{cases}
		\end{equation}
		
		\item[(ii)] For any $(n,m)\in\bbN^2$, the type II multiple orthogonal polynomials are related via
		\begin{equation}\label{eq:TypeII-XM1}
			P_{n,m}(x) = 
			\begin{cases}
				\widehat{P}_{n,m}(x), & \mbox{if } n \ge m, \\
				\widetilde{P}_{m,n}(x), & \mbox{if } n \le m.
			\end{cases}
		\end{equation}
		
		\item[(iii)]  For any $(n,m)\in\bbN^2$, the type I multiple orthogonal polynomials are related via
		\begin{equation}\label{eq:TypeI-XM1-Right}
			\begin{aligned}
				A_{n,m}(x) &= \widehat{A}_{n,m}(x) + c \widehat{B}_{n,m}(x), \\
				B_{n,m}(x) &= -\widehat{B}_{n,m}(x),
			\end{aligned}
			\qquad \text{if } n \ge m,
		\end{equation}
		and
		\begin{equation}\label{eq:TypeI-XM1-Left}
			\begin{aligned}
				A_{n,m}(x) &= \widetilde{B}_{m,n}(x), \\
				B_{n,m}(x) &= \widetilde{A}_{m,n}(x) - \frac{1}{c} \widetilde{B}_{m,n}(x),
			\end{aligned}
			\qquad \text{if } n \le m,
		\end{equation}
		where $c = \int_{\bbR_-} d\sigma_2$.
	\end{itemize}
\end{theorem}


\begin{proof}
	\textit{(i)} Applying the Andreief identity to~\eqref{eq:Hankel}, we obtain 
	\begin{equation}\label{eq:Hproof1}
		H_{n,m}[\bm\mu]
		=
		\frac{1}{N!} \int_{\bm{x}\in\bbR_+^N}   
		U_{n,m}(\bm{x})
		\prod_{1 \le j<k \le N}(x_k-x_j) \, \prod_{k=1}^N d\sigma_1(x_k),
	\end{equation}
	where $U_{n,m}(\bm{x}) = U_{\bm{x}} \big(1,x,\ldots,x^{n-1}, x \mathfrak{m}_{\sigma_2}(x), \ldots,x^{m} \mathfrak{m}_{\sigma_2}(x)\big)$. Because the integration measure and the Vandermonde determinant $\Delta(\bm{x}) = \prod_{j<k}(x_k-x_j)$ are  positive on the ordered simplex $x_1 < \ldots < x_N$, it suffices to show that $U_{n,m}(\bm{x})$ maintains a constant, non-zero sign on this domain.
	
	\textit{Case $n \ge m$}.
	Recall $d\widehat{\sigma}_2(t) = -t \, d\sigma_2(t)$. Using Lemma \ref{lem:transforms}(i) with $z_0=0$, we can expand the functions in the bottom block of  $U_{n,m}(\bm{x})$ as
	\begin{equation}
		x^j \mathfrak{m}_{\sigma_2}(x) = x^{j-1} (x\mathfrak{m}_{\sigma_2}(x)) = -x^{j-1} \mathfrak{m}_{\widehat{\sigma}_2}(x) + x^{j-1} \int_{\bbR_-} d\sigma_2(t)
	\end{equation}
	for each $j=1, \dots, m$. Because $j-1 \le m-1 \le n-1$, the polynomial terms $x^{j-1} \int d\sigma_2$ are spanned by the top block of $n$ rows. Applying elementary row operations, we reduce $U_{n,m}(\bm{x})$ to $U_{\bm{x}}(1,\ldots,x^{n-1},-\mathfrak{m}_{\widehat{\sigma}_2},\ldots,-x^{m-1} \mathfrak{m}_{\widehat{\sigma}_2})$. 
	
	Factoring out the minus sign from each of the bottom $m$ rows yields $(-1)^m U_{\bm{x}}[\mathcal{N}(\widehat{\sigma}_1, \widehat{\sigma}_2)]$, where $U_{\bm{x}}[\mathcal{N}]$ is the corresponding $U_{\bm{x}}$ determinant for the 
	system 
	$\widehat{\bm\mu}$ at the index $(n,m)$. This proves $H_{n,m} = 
		(-1)^{m} \widehat{H}_{n,m}$. 


	\textit{Case $n \le m$}.
	We use~\eqref{eq:Hproof1}. 
	 We extract a factor $x_k \mathfrak{m}_{\sigma_2}(x_k)$ from each of its columns to write
	\begin{multline}\label{eq:Uxproof}
		H_{n,m}[\bm\mu]
		=
		\frac{1}{N!} \int_{\bm{x}\in\bbR_+^N}   
		\Big(\prod_{k=1}^N x_k \mathfrak{m}_{\sigma_2}(x_k)\Big) U_{\bm{x}}\Big(\frac{1}{x\mathfrak{m}_{{\sigma}_2}}, \ldots, \frac{x^{n-1}}{x\mathfrak{m}_{{\sigma}_2}}, 1, \ldots, x^{m-1}\Big)
		\\
		\times
		\prod_{1 \le j<k \le N}(x_k-x_j) \, \prod_{k=1}^N d\sigma_1(x_k),
	\end{multline}
	
	Now recall~\eqref{eq:sigma2tilde}:
	\begin{equation}\label{eq:1/xm}
		\frac{1}{x \mathfrak{m}_{\sigma_2}(x) } = \frac{1}{c} + \mathfrak{m}_{\widetilde{\sigma}_2}(x).
	\end{equation}
	Using this relation and performing 
	 elementary row operations (utilizing $j \le n-1 \le m-1$) in the $U_{\bm{x}}$ determinant in~\eqref{eq:Uxproof}, we can reduce it  to $U_{\bm{x}}(\mathfrak{m}_{\widetilde{\sigma}_2}, x \mathfrak{m}_{\widetilde{\sigma}_2}, \ldots, x^{n-1} \mathfrak{m}_{\widetilde{\sigma}_2}, 1, \ldots, x^{m-1})$.

	Now combine $x_k \mathfrak{m}_{\sigma_2}(x_k) d\sigma_1(x_k)$ into $d\widetilde{\sigma}_1(x_k)$ and swap the top block (size $n$) of the $U_{\bm{x}}$ determinant and its bottom block (size $m$). This introduces a sign of $(-1)^{nm}$ and proves $H_{n,m} = (-1)^{nm} \widetilde{H}_{m,n}$.
	

	\smallskip
	
	\smallskip
	
	\textit{(ii)} \textit{Case $n \ge m$}. 
	By definition, $P_{n,m}(x)$ is a monic polynomial of degree $n+m$ satisfying orthogonality against $1,\ldots,x^{n-1}$ with respect to $\sigma_1=\widehat{\sigma}_1$, and against $1,\ldots,x^{m-1}$ with respect to $x\mathfrak{m}_{\sigma_2}(x){\sigma}_1 = x\mathfrak{m}_{\sigma_2}(x)\widehat{\sigma}_1$.
	
	Because $\widehat{\sigma}_2 = -t\sigma_2$, we have $x\mathfrak{m}_{\sigma_2}(x) = c - \mathfrak{m}_{\widehat{\sigma}_2}(x)$ by Lemma~\ref{lem:transforms}\textit{(i)}. Substituting this into the conditions for the second weight yields:
	\begin{equation}
		\int_{\bbR} P_{n,m}(x) x^j \big(c - \mathfrak{m}_{\widehat{\sigma}_2}(x)\big) d\widehat{\sigma}_1(x) = 0, \qquad j=0, \dots, m-1.
	\end{equation}
	Because $j \le m-1 \le n-1$, the integral $c \int P_{n,m}(x) x^j d\widehat{\sigma}_1(x)$ is zero due to the orthogonality conditions on the first weight. This implies $P_{n,m}$ satisfies all the type II orthogonality for  the system $\widehat{\bm{\mu}}$. Since it is perfect, we have $P_{n,m}=\widehat{P}_{n,m}$ by uniqueness. 
	
	\textit{Case $n \le m$}. 
	Recall $\widetilde{\sigma}_1 = x \mathfrak{m}_{\sigma_2}(x) \sigma_1$. Substituting 
	\begin{align*}
	& d\sigma_1(x)  = (x\mathfrak{m}_{\sigma_2}(x))^{-1} d\widetilde{\sigma}_1(x) = \big(c^{-1} + \mathfrak{m}_{\widetilde{\sigma}_2}(x)\big) d\widetilde{\sigma}_1(x),
	\\
	 & x \m_{\sigma_2}(x) d\sigma_1(x)  =  d\widetilde{\sigma}_1(x),
	\end{align*}
	 into the orthogonality conditions for $\bm{\mu}$ yields:
	\begin{alignat}{3}
		& \int_{\bbR} P_{n,m}(x) x^k \big(c^{-1} + \mathfrak{m}_{\widetilde{\sigma}_2}(x)\big) d\widetilde{\sigma}_1(x) = 0, \qquad && k=0, \dots, n-1, \\
		& \int_{\bbR} P_{n,m}(x) x^j d\widetilde{\sigma}_1(x) = 0, \qquad && j=0, \dots, m-1.
	\end{alignat}
	Because $k \le n-1 \le m-1$, the integral $c^{-1} \int P_{n,m}(x) x^k d\widetilde{\sigma}_1(x)$ vanishes via the second set of conditions. This shows $P_{n,m} = \widetilde{P}_{m,n}$ by perfectness of $\widetilde{\bm\mu}$.

	\smallskip
	
	\textit{(iii)} \textit{Case $n \ge m$}. 
	The type I vector $(A_{n,m}, B_{n,m})$ for $\bm{\mu}$ satisfies $\deg A_{n,m} \le n-1$, $\deg B_{n,m} \le m-1$, and
	\begin{equation}\label{eq:proofTypI}
		\int_{\bbR} \big(A_{n,m}(x) + B_{n,m}(x) x\mathfrak{m}_{\sigma_2}(x)\big) x^k d\sigma_1(x) = \delta_{k, n+m-1}, \quad k=0, \dots, n+m-1.
	\end{equation}
	Substituting $x\mathfrak{m}_{\sigma_2}(x) = c - \mathfrak{m}_{\widehat{\sigma}_2}(x)$ and rearranging the terms gives:
	\begin{equation}
		\int_{\bbR} \Big( \big[A_{n,m}(x) + cB_{n,m}(x)\big] - B_{n,m}(x) \mathfrak{m}_{\widehat{\sigma}_2}(x) \Big) x^k d\widehat{\sigma}_1(x) = \delta_{k, n+m-1}.
	\end{equation}
	Because $n \ge m$, the polynomial $A_{n,m} + cB_{n,m}$ has degree at most $n-1$, and $-B_{n,m}$ has degree at most $m-1$. These match the degree requirements for the type I polynomial $(\widehat{A}_{n,m}, \widehat{B}_{n,m})$ of $\widehat{\bm{\mu}}$. By normality, $\widehat{A}_{n,m} = A_{n,m} + cB_{n,m}$ and $\widehat{B}_{n,m} = -B_{n,m}$, which inverts to~\eqref{eq:TypeI-XM1-Right}.
	
	\textit{Case $n \le m$}. 
	Factoring out $x\mathfrak{m}_{\sigma_2}(x)$ from~\eqref{eq:proofTypI} and using~\eqref{eq:1/xm} yields:
	\begin{equation}
		\int_{\bbR} \big(A_{n,m}(x) [c^{-1} + \mathfrak{m}_{\widetilde{\sigma}_2}(x)] + B_{n,m}(x)\big) x^k d\widetilde{\sigma}_1(x) = \delta_{k, n+m-1}.
	\end{equation}
	Rearranging the integrand gives:
	\begin{equation}
		\int_{\bbR} \Big( \big[B_{n,m}(x) + c^{-1}A_{n,m}(x)\big] + A_{n,m}(x)\mathfrak{m}_{\widetilde{\sigma}_2}(x) \Big) x^k d\widetilde{\sigma}_1(x) = \delta_{k, n+m-1}.
	\end{equation}
	Because $n \le m$, the polynomial $B_{n,m} + c^{-1}A_{n,m}$ has degree at most $m-1$, and $A_{n,m}$ has degree at most $n-1$. This implies $\widetilde{A}_{m,n} = B_{n,m} + c^{-1}A_{n,m}$ and $\widetilde{B}_{m,n} = A_{n,m}$, which inverts to~\eqref{eq:TypeI-XM1-Left}.
\end{proof}

\begin{proof}[Proof of Theorem \ref{thm:1-xm_Intro}]
	By translation and reflection, we can reduce the setting of the theorem to the case $I_1 \subseteq \bbR_+$, $I_2 \subseteq \bbR_-$, and $x_0=0$. 
	
	By Theorem \ref{thm:XM1_Properties}\textit{(i)}, the Hankel determinant $H_{n,m}[\bm{\mu}]$ is equal (up to a non-zero sign) to the Hankel determinant of a standard Nikishin system $\widehat{\bm\mu}$ or $\widetilde{\bm\mu}$. Because those are perfect, $H_{n,m}[\bm{\mu}] \neq 0$ for all $(n,m) \in \bbN^2$.
\end{proof}

\begin{remark}\label{rem:NNRR1}
	The identities established in Theorem~\ref{thm:XM1_Properties} then imply that the recurrence coefficients of $\bm{\mu}$ coincide with those of $\widehat{\bm{\mu}}$ and $\widetilde{\bm{\mu}}$ in the respective cones. 
	Let $p_{n,m}$ denote the subleading coefficient of $P_{n,m}(x) = x^{n+m} + p_{n,m} x^{n+m-1} + \dots$. By expanding~\eqref{eq:nnrr-1}--\eqref{eq:nnrr-2} relations, one obtains $b_{n,m;1} = p_{n,m} - p_{n+1,m}$ and $b_{n,m;2} = p_{n,m} - p_{n,m+1}$. 
	Because $P_{n,m}(x) = \widehat{P}_{n,m}(x)$ for $n \ge m$, we have $p_{n,m} = \widehat{p}_{n,m}$ and therefore
	\begin{equation*}
		b_{n,m;1} = \widehat{b}_{n,m;1} \quad (n \ge m), \qquad b_{n,m;2} = \widehat{b}_{n,m;2} \quad (n \ge m+1).
	\end{equation*}
	Applying $P_{n,m}(x) = \widetilde{P}_{m,n}(x)$ for $n \le m$, we get
	\begin{equation*}
		b_{n,m;1} = \widetilde{b}_{m,n;2} \quad (n \le m-1), \qquad b_{n,m;2} = \widetilde{b}_{m,n;1} \quad (n \le m).
	\end{equation*}
	For the $a$-coefficients, we use~\eqref{eq:a-coeffs} which imply
	\begin{align*}
		a_{n,m;1} &= \widehat{a}_{n,m;1}, \qquad a_{n,m;2} = \widehat{a}_{n,m;2}, \qquad &&\text{for } n \ge m+1, \\
		a_{n,m;1} &= \widetilde{a}_{m,n;2}, \qquad a_{n,m;2} = \widetilde{a}_{m,n;1}, \qquad &&\text{for } n \le m-1.
	\end{align*}
	Remarkably, every $b$-coefficient of $\bm\mu$ is inherited from one of the two unperturbed Nikishin systems $\widehat{\bm\mu}$, $\widetilde{\bm\mu}$. The same holds true for the $a$-coefficients, with the exception of the  diagonal coefficients $n=m$ which cross-couple the two systems.
\end{remark}

\begin{corollary}\label{cor:TypeI-RealRoots}
	In the setting of Theorem~\ref{thm:XM1_Properties}:
	\begin{itemize}
		\item[(i)] For any $(n,m)\in\bbN^2$, the type II polynomial $P_{n,m}(x)$ is real-rooted. Its $n+m$ simple zeros belong to $\mathring{I}_1$ and interlace with the zeros of each $P_{n+1,m}(x)$ and $P_{n,m+1}(x)$. 
		
		\item[(ii)] For $n+1 \le m$, $n>0$, the type I polynomial $A_{n,m}(x)$ is real-rooted. Its $n-1$ simple zeros belong to the interior of $\operatorname{conv}(I_2 \cup \{0\})$ and interlace with the zeros of each $A_{n-1,m}(x)$ and $A_{n,m-1}(x)$.
		
		\item[(iii)] For $m+1 \le n$, $m>0$, the type I polynomial $B_{n,m}(x)$ is real-rooted. Its $m-1$ simple zeros belong to $\mathring{I}_2$ and interlace with the zeros of each $B_{n-1,m}(x)$ and $B_{n,m-1}(x)$.
	\end{itemize}
\end{corollary}

\begin{remark}\label{rem:TypeI-Roots}
	The results of this corollary can be improved further. This falls outside the scope of the present work and will be the subject of a separate paper.
\end{remark}

\begin{proof}
	\textit{(i)} At every location $(n,m)$ the system is AT on $I_1$. This immediately implies the statement, see~\cite{HanVA,Ismail}.
	
	\textit{(ii)} For $n+1 \le m$, Theorem~\ref{thm:XM1_Properties}\textit{(iii)} gives $A_{n,m} = \widetilde{B}_{m,n}$. Since $\widetilde{\bm{\mu}} = \mathcal{N}(\widetilde{\sigma}_1, \widetilde{\sigma}_2)$ is a Nikishin system with $\supp\,\widetilde{\sigma}_2 \subseteq \operatorname{conv}(I_2 \cup \{0\})$ and $n+1 \le m$, both the zero location and the stated interlacing relations follow from~\cite{KVInterlacing}.
	
	\textit{(iii)} For $m+1 \le n$ Theorem~\ref{thm:XM1_Properties}\textit{(iii)} gives $B_{n,m} = -\widehat{B}_{n,m}$. Since $\widehat{\bm{\mu}} = \mathcal{N}(\widehat{\sigma}_1, \widehat{\sigma}_2)$ is a Nikishin system with $\operatorname{conv}(\supp\,\widehat{\sigma}_2) = I_2$ and $m+1 \le n$, both the zero location and the stated interlacing relations follow from~\cite{KVInterlacing}.
\end{proof}

\subsection{System $(x \sigma_1, \mathfrak{m}_{\sigma_2}\sigma_1)$}


\begin{proof}[Proof of Theorem \ref{thm:2-xm_Intro}]
	By translation and reflection, we can reduce the setting of the theorem to the case $\supp\,\sigma_1 \subseteq \bbR_+$, $\supp\,\sigma_2 \subseteq \bbR_-$, and $x_0=0$. Therefore, we consider the modified system $\bm{\nu} = (x\sigma_1, \mathfrak{m}_{\sigma_2} \sigma_1)$.

\textit{Case $n\ge m-2$.} If $m=0$, then $H_{n,0}[\bm{\nu}]$ is the standard Hankel determinant of the positive measure $x\sigma_1$ on $\bbR_+$, so $H_{n,0}[\bm{\nu}] > 0$. Assume now that $m \ge 1$ and let $N = n+m$. 
Applying the Andr\'eief identity to~\eqref{eq:Hankel}, we obtain the integral representation~\cite{Joh,Kui}:
\begin{multline}\label{eq:Htemp2}
	H_{n,m}[\bm{\nu}]
	\\=
	\frac{1}{N!} \int_{\bm{x}\in\bbR_+^N}   
	U_{\bm{x}} \Big(x,x^2,\ldots,x^{n},
	\mathfrak{m}_{\sigma_2}(x),x \mathfrak{m}_{\sigma_2}(x),\ldots,x^{m-1} \mathfrak{m}_{\sigma_2}(x)\Big)
	\\
	\times \prod_{j<k}(x_k-x_j) \, \prod_{j=1}^N d\sigma_1(x_j).
\end{multline}

Assume for now that $\int_{\bbR_-} |t|^{-1} d\sigma_2(t) < \infty$. Since $\frac{x^j}{x-t} = \frac{x t^{j-1}}{x-t} + \sum_{i=1}^{j-1} t^{j-1-i} x^i$ and $m-2\le n$, by performing elementary row operations using the first $n$ rows for each $j=1,\ldots,m-1$, we can reduce the row with entries $x_k^j \mathfrak{m}_{\sigma_2}(x_k)$ in the determinant to $\int_{\bbR_-} \frac{ x_k t^{j-1} d\sigma_2(t)}{x_k-t}$. Extracting $x_k$ from each column we obtain
	\begin{multline}\label{eq:Htemp}
		H_{n,m}[\bm{\nu}]
		\\=
		\frac{1}{N!} \int_{\bm{x}\in\bbR_+^N}   
		U_{\bm{x}} \Big(1,x,\ldots,x^{n-1},
		\int \frac{d\sigma_2(t)}{x(x-t)} ,\int \frac{ t^{0} d\sigma_2(t)}{x-t} ,\ldots,\int \frac{ t^{m-2} d\sigma_2(t)}{x-t} \Big)
		\\
		\times \prod_{j<k}(x_k-x_j) \, \prod_{k=1}^N x_k \, \prod_{j=1}^N d\sigma_1(x_j).
	\end{multline}
	Using the identity $\frac{1}{x(x-t)} =  \frac{1}{t(x-t)} - \frac{1}{tx}$, we can split the last $U_{\bm{x}}$ determinant as $U_1 - U_2$, where
	\begin{align}
		U_1 & = U_{\bm{x}} \Big(1,x,\ldots,x^{n-1},
		\int \frac{t^{-1}d\sigma_2(t)}{x-t} ,\int \frac{ t^{0} d\sigma_2(t)}{x-t},\ldots,\int \frac{ t^{m-2} d\sigma_2(t)}{x-t} \Big),
		\\
		U_2& =    U_{\bm{x}} \Big(1,x,\ldots,x^{n-1},
		\int \frac{d\sigma_2(t)}{xt} ,\int \frac{ t^{0} d\sigma_2(t)}{x-t},\ldots,\int \frac{ t^{m-2} d\sigma_2(t)}{x-t} \Big).
	\end{align}
	
	Plugging $U_{\bm{x}} = U_1 - U_2$ into~\eqref{eq:Htemp}, we obtain two integrals, the first of which reduces to
	\begin{multline}\label{eq:u1term}
		\frac{1}{N!} \int_{\bm{x}\in\bbR_+^N}   
		U_{1} \prod_{j<k}(x_k-x_j) \, \prod_{k=1}^N x_k \, \prod_{j=1}^N d\sigma_1(x_j)
		\\
		=
		\frac{1}{N! m!} \int_{\bm{t}\in\bbR_-^m}    \int_{\bm{x}\in\bbR_+^N}   
		U_{\bm{x}} \left( 1,x,\ldots,x^{n-1}, \frac{1}{x-t_1},\ldots,\frac{1}{x-t_m}\right) 
		\\
		\times \prod_{j=1}^m t^{-1}_j \, \prod_{j<k}(t_k-t_j) \, \prod_{j<k}(x_k-x_j) \, \prod_{k=1}^N x_k \, \prod_{j=1}^N d\sigma_1(x_j) \, \prod_{j=1}^m d\sigma_2(t_j)
	\end{multline}
	by the generalized Andréief identity (see, e.g.,~\cite[Prop. 3.1]{KVMLOPUC1}).
	
	For the integral with $U_2$, we first extract the integral $\int t_1^{-1} d\sigma_2(t_1)$ from the $(n+1)$-st row. We then apply the generalized Andréief identity to the remaining $m-1$ integral rows, yielding 
	\begin{multline}\label{eq:u2term}
		\frac{1}{N!} \int_{\bm{x}\in\bbR_+^N}   
		U_{2} \prod_{j<k}(x_k-x_j) \, \prod_{k=1}^N x_k \, \prod_{j=1}^N d\sigma_1(x_j)
		\\
		=\frac{1}{N! (m-1)!} \int_{\bm{t}\in\bbR_-^m}   \int_{\bm{x}\in\bbR_+^N}   
		U_{\bm{x}} \left( 1,x,\ldots,x^{n-1}, \frac{1}{x},\frac{1}{x-t_2}\ldots,\frac{1}{x-t_m}\right) 
		\\
		\times t_1^{-1} \, \prod_{2\le j<k\le m}(t_k - t_j) 
		\prod_{j<k}(x_k-x_j) \, \prod_{k=1}^N x_k \, \prod_{j=1}^N d\sigma_1(x_j) \, \prod_{j=1}^m d\sigma_2(t_j).
	\end{multline}
	Observe that the $U_{\bm{x}}$ determinant in~\eqref{eq:u1term} is exactly $C(\bm{x};\bm{t})$ from~\eqref{eq:C}, while the $U_{\bm{x}}$ determinant in~\eqref{eq:u2term} is  $C(\bm{x};0,t_2,\ldots,t_m)$. 
	
	Note that the second integral is not symmetric in the $(t_1,\ldots,t_m)$ variables, so we now perform $\bm{t}$-symmetrization by replacing its integrand by its average over the permutation group $S_m$ applied to the $\bm{t}$ variables. Let us introduce the notation $\bm{t}_{\setminus r}$ for the ordered set $(t_1,\ldots,t_m)\setminus\{t_r\}$ and $V(\bm{v})$ for the Vandermonde determinant $\prod_{j<k} (v_k-v_j)$, and $V(\bm{v};\bm{u}) = \prod \prod (v_j - u_k)$.     
	
	Now group the sum over $S_m$ by the index $r = \sigma(1)$ that occupies the $t_1$ position. For a fixed $r$, the remaining variables undergo permutations in $S_{m-1}$. Because the product $t_r^{-1} C(\bm{x}; 0, \bm{t}_{\setminus r})V( \bm{t}_{\setminus r}) $ is  symmetric with respect to those $m-1$ variables, the sum over $S_{m-1}$ simply produces a factor of $(m-1)!$. Therefore the symmetrization of ~\eqref{eq:u2term} is
	\begin{multline}\label{eq:u2term2}
		\frac{1}{N!m!} \int_{\bm{t}\in\bbR_-^m}   \int_{\bm{x}\in\bbR_+^N}   
		\left( \sum_{r=1}^m  \tfrac{1}{t_r} C(\bm{x}; 0, \bm{t}_{\setminus r})V( \bm{t}_{\setminus r}) \right) 
		\\
		\times 
		\prod_{j<k}(x_k-x_j) \, \prod_{k=1}^N x_k \, \prod_{j=1}^N d\sigma_1(x_j) \, \prod_{j=1}^m d\sigma_2(t_j).
	\end{multline}
	
	Combining~\eqref{eq:u1term} and~\eqref{eq:u2term2} we obtain
	\begin{multline}
		\frac{1}{N!m!} \int_{\bm{t}\in\bbR_-^m}   \int_{\bm{x}\in\bbR_+^N}   
		\left( \Big( \prod_{k=1}^m t_k^{-1} \Big) C(\bm{x}; \bm{t})V( \bm{t}) -  \sum_{r=1}^m  t_r^{-1} C(\bm{x}; 0, \bm{t}_{\setminus r})V( \bm{t}_{\setminus r}) \right) 
		\\
		\times 
		\prod_{j<k}(x_k-x_j) \, \prod_{k=1}^N x_k \, \prod_{j=1}^N d\sigma_1(x_j) \, \prod_{j=1}^m d\sigma_2(t_j).
	\end{multline}
	
	Substitute the Cauchy-Vandermonde formulas~\eqref{eq:C}: $C(\bm{x};\bm{t}) = (-1)^{S(n,m)} \frac{V(\bm{x}) V(\bm{t})}{V(\bm{x};\bm{t})}$ and $C(\bm{x};0, \bm{t}_{\setminus r}) = (-1)^{S(n,m)} \frac{V(\bm{x}) V(\bm{t}_{\setminus r}) \prod_{i \neq r} t_i}{V(\bm{x};\bm{t}_{\setminus r}) \prod x_j}$. 
	
	Using the identities $V(\bm{t}_{\setminus r})^2 = V(\bm{t})^2 / \prod_{j \neq r}(t_r - t_j)^2$ and $V(\bm{x};\bm{t}_{\setminus r}) = V(\bm{x};\bm{t})/ \prod_{j=1}^N (x_j-t_r)$, we finally arrive at the integrand
	\begin{multline}\label{eq:jpdf}
		\frac{(-1)^{S(n,m)}}{N!m!}   
		\left( 
		1
		-  \sum_{r=1}^m   \prod_{j=1}^N \Big(1-\frac{t_r}{x_j} \Big)
		\prod_{k\ne r} \Big(\frac{t_k}{t_r - t_k} \Big)^2 
		\right) 
		\\
		\times 
		\frac{V(\bm{t})^2 V(\bm{x})^2}{V(\bm{x};\bm{t})} \,   \prod_{j=1}^N x_j d\sigma_1(x_j) \, \prod_{j=1}^m t_j^{-1} d\sigma_2(t_j).
	\end{multline}
	
	Because the distribution is symmetric we can remove the $N! m!$ normalization and restrict the domain to a chamber $t_1<t_2<\ldots<t_m<0$ and $0<x_1<\ldots<x_N$.
	
	Since $t_k<t_m<0$ for any $k < m$ and $x_j > 0$, we see that 
	\begin{equation}
		\frac{t_k^2}{(t_m - t_k)^2} > 1, \qquad 1-\frac{t_m}{x_j}>1.
	\end{equation}
	Because every individual factor in the $r=m$ term of the sum in~\eqref{eq:jpdf} is $>1$, the product itself is $>1$. Since all other terms in the sum are positive, the expression in the brackets is strictly negative. Taking into account the sign $(-1)^m$ from $\prod_{j=1}^m t_j^{-1}$, we obtain for all $m \ge 1$ and $n \ge m-2$,
	\begin{equation}
		\operatorname{sgn}(H_{n,m}[\bm{\nu}]) = (-1)^{S(n,m)+m+1} \ne 0.
	\end{equation}
	
	Finally, if $\int |t|^{-1}d\sigma_2(t)$ diverges, we define the truncated measure $d\sigma_{2,\epsilon} = \chi_{(-\infty, -\epsilon]} d\sigma_2$. For each $\epsilon > 0$, the above result holds true. By the monotone convergence theorem applied to the expression \eqref{eq:jpdf}, the magnitude of the integral increases as $\epsilon \to 0$. Since it is non-zero for any $\epsilon > 0$, the limit is bounded away from zero. 
	
	\smallskip
	
	\textit{Case $n\le m-2$}.  We use~\eqref{eq:Htemp2} and extract $\mathfrak{m}_{\sigma_2}(x_k)$ from each column of the matrix to write 
	\begin{multline}\label{eq:Htemp3}
		H_{n,m}[\bm{\nu}]
		\\=
		\frac{1}{N!} \int_{\bm{x}\in\bbR_+^N}   
		U_{\bm{x}} \Big(\frac{x}{\mathfrak{m}_{\sigma_2}(x)},\frac{x^2}{\mathfrak{m}_{\sigma_2}(x)},\ldots,\frac{x^{n}}{\mathfrak{m}_{\sigma_2}(x)},
		1,x ,\ldots,x^{m-1} \Big)
		\\
		\times \prod_{j<k}(x_k-x_j) \, \prod_{k=1}^N \mathfrak{m}_{\sigma_2}(x_k) d\sigma_1(x_k).
	\end{multline}
	
	Swapping the top $n$-row block with the bottom $m$-row block results in a factor of $(-1)^{nm}$. Now 
	use the inversion $\frac{x^j}{\mathfrak{m}_{\sigma_2}(x)} = b x^j + cx^{j+1} -x^j \mathfrak{m}_{\sigma_2^{(1)}}(x)$ (Lemma~\ref{lem:transforms}\textit{(iii)}) for $j=1,\ldots,n$. Because $n \le m-2$, the maximal degree $n+1$ is bounded by $m-1$. Thus, we can perform elementary row operations using the bottom block to eliminate the $b$- and $c$-terms. Extracting the remaining minus sign gives 
	\begin{multline}\label{eq:Htemp4}
		H_{n,m}[\bm{\nu}]
		\\=
		\frac{(-1)^{n(m+1)}}{N!} \int_{\bm{x}\in\bbR_+^N}   
		U_{\bm{x}} \Big(1,x ,\ldots,x^{m-1}, x\mathfrak{m}_{\sigma^{(1)}_2}, x^2\mathfrak{m}_{\sigma^{(1)}_2},\ldots, x^n\mathfrak{m}_{\sigma^{(1)}_2} \Big)
		\\
		\times \prod_{j<k}(x_k-x_j) \, \prod_{k=1}^N \mathfrak{m}_{\sigma_2}(x_k) d\sigma_1(x_k),
	\end{multline}
	which is exactly $(-1)^{n(m+1)} H_{m,n}[(\mathfrak{m}_{\sigma_2} \sigma_1, x\mathfrak{m}_{\sigma_2^{(1)}} \mathfrak{m}_{\sigma_2} \sigma_1)]$. 
	By Theorem~\ref{thm:XM1_Properties}, the modified Nikishin system $(\mathfrak{m}_{\sigma_2} \sigma_1, x\mathfrak{m}_{\sigma_2^{(1)}} \mathfrak{m}_{\sigma_2} \sigma_1)$ is perfect, which completes the proof.
\end{proof}

\begin{remark}\label{rem:RMT}
	The symmetrized integrand in~\eqref{eq:jpdf} is sign-definite, and therefore dividing it by $H_{n,m}[\bm{\nu}]$ yields a valid joint probability density function (jpdf). In random matrix theory, this distribution corresponds to a multiple orthogonal polynomial ensemble (MOPE), interpreted as a two-component interacting log-gas with $N$ particles at $\bm{x} \in \bbR_+^N$ and $m$ particles at $\bm{t} \in \bbR_-^m$ \cite{Bor, Kui}. 
	
	Moreover, this distribution can be viewed as a perturbation of the classical Nikishin MOPE associated with $\mathcal{N}(x\sigma_1(x), |t|^{-1}\sigma_2(t))$. If the additive perturbation term (the bracketed sum in~\eqref{eq:jpdf}) is removed, the density reduces exactly to the jpdf of this unperturbed Nikishin ensemble (cf.~\cite{Kui}).
\end{remark}

\begin{remark}\label{rem:NNRR2}
	Combine the reduction in~\eqref{eq:Htemp4} for $n \le m-2$ with Theorem~\ref{thm:XM1_Properties} applied to the system $(\m_{\sigma_2}\sigma_1, x\m_{\sigma_2^{(1)}}\m_{\sigma_2}\sigma_1)$ at index $(m,n)$ (where $m \ge n+2 > n$). This yields $P_{n,m}(x) = \breve{P}_{m,n}(x)$ and $H_{n,m}[\bm{\nu}] = (-1)^n \breve{H}_{m,n}$ for $n \le m-2$, where $\breve{\bm{\nu}} = \mathcal{N}(\m_{\sigma_2}\sigma_1, -t\sigma_2^{(1)})$. Analogously to Remark~\ref{rem:NNRR1}, this allows us to extract the recurrence coefficients in this cone:
	\begin{equation*}
		b_{n,m;2} = \breve{b}_{m,n;1} \quad (n \le m-2), \qquad b_{n,m;1} = \breve{b}_{m,n;2} \quad (n \le m-3),
	\end{equation*}
	and
	\begin{equation*}
		a_{n,m;1} = \breve{a}_{m,n;2}, \qquad a_{n,m;2} = \breve{a}_{m,n;1}, \qquad \text{for } n \le m-3.
	\end{equation*}
\end{remark}

\begin{corollary}\label{cor:TypeI-RealRoots2}
	In the setting of Theorem~\ref{thm:2-xm_Intro}:
	\begin{itemize}
		\item[(i)] For any $(n,m)\in\bbN^2$, the type II polynomial $P_{n,m}(x)$ is real-rooted. Its $n+m$ simple zeros belong to $\mathring{I}_1$ and interlace with the zeros of each $P_{n+1,m}(x)$ and $P_{n,m+1}(x)$. 
		
		\item[(ii)] For $n \le m-3$, $n>0$, the type I polynomial $A_{n,m}(x)$ is real-rooted. Its $n-1$ simple zeros belong to $\mathring{I}_2$ and interlace with the zeros of each $A_{n-1,m}(x)$ and $A_{n,m-1}(x)$.
	\end{itemize}
\end{corollary}

\begin{remark}\label{rem:TypeI-Roots2}
	Determining the zero location and interlacing properties for $B_{n,m}(x)$ cannot be achieved through this mapping since when $n\ge m-1$ the system is only AT and not Nikishin-reducible. 
\end{remark}

\begin{proof}
	\textit{(i)} The proof of Theorem~\ref{thm:2-xm_Intro} shows that $U_{\bm{x}}\ne 0$ for every strictly increasing $\bm{x} \subset \mathring{I}_1$. This implies that the system is not only perfect, but is an AT system on $I_1$ at every $(n,m)\in\bbN^2$. Property \textit{(i)} therefore follows immediately from standard AT theory (see, e.g.,~\cite{HanVA,Ismail}).
	
	\textit{(ii)} Let $n \le m-2$. The type I orthogonality conditions for $\bm{\nu}$ yield
	\begin{multline}
		\int_{I_1} \Big( A_{n,m}(x) \frac{x}{\mathfrak{m}_{\sigma_2}(x)} + B_{n,m}(x) \Big) x^k \mathfrak{m}_{\sigma_2}(x) d\sigma_1(x) 
		\\
		= \delta_{k,n+m-1}, \quad k=0, \dots, n+m-1.
	\end{multline}
	Using Lemma~\ref{lem:transforms}\textit{(i)} and \textit{(iii)}, we have $\frac{x}{\mathfrak{m}_{\sigma_2}(x)} = cx^2 + bx - x\mathfrak{m}_{\sigma_2^{(1)}}(x) = cx^2 + bx - d + \mathfrak{m}_{-t\sigma_2^{(1)}}(x)$, where $c = (\int d\sigma_2)^{-1}$ and $d = \int d\sigma_2^{(1)}$. Defining $d\breve{\nu}_1 = \mathfrak{m}_{\sigma_2} d\sigma_1$, we can rewrite the orthogonality condition as
	\begin{multline}\label{eq:TypeI-Mapped}
		\int_{I_1} \Big( \underbrace{\big[B_{n,m}(x) + A_{n,m}(x)(cx^2+bx-d)\big]}_{\breve{A}_{m,n}(x)} + \underbrace{A_{n,m}(x)}_{\breve{B}_{m,n}(x)} \mathfrak{m}_{-t\sigma_2^{(1)}}(x) \Big) x^k d\breve{\nu}_1(x) 
		\\
		= \delta_{k,n+m-1}, \quad k=0, \dots, n+m-1.
	\end{multline}
	Because $\deg A_{n,m} \le n-1$ and $n \le m-2$, the polynomial $\breve{A}_{m,n}$ has degree bounded by $\max(m-1, n+1) = m-1$. 
	Thus,~\eqref{eq:TypeI-Mapped} is the type I orthogonality condition for the standard Nikishin system $\breve{\bm{\nu}} = \mathcal{N}(\breve{\nu}_1, -t\sigma_2^{(1)})$ at the swapped multi-index $(m, n)$. 
	For $n \le m-3$, the indices $(m,n)$, $(m,n-1)$, and $(m-1,n)$ all satisfy this reduction, so~\cite{KVInterlacing} guarantees that $\breve{B}_{m,n}(x) = A_{n,m}(x)$ has $n-1$ simple zeros in $\mathring{I}_2$ that interlace with those of $A_{n-1,m}(x)$ and $A_{n,m-1}(x)$.
\end{proof}


\section{Signs of Hankel determinants and recurrence coefficients}\label{ss:Sign}

In this section we determine the signs of the stacked Hankel determinants and the nearest-neighbor recurrence coefficients for a Nikishin system and both of its partial Christoffel modifications.

\begin{proposition}\label{prop:Hsigns}
	Suppose that $\supp\, \sigma_1 \subseteq \bbR_+$ and $\supp\, \sigma_2 \subseteq \bbR_-$.
	\begin{itemize}
		\item[(i)]  For the Nikishin system $\mathcal{N}(\sigma_1,\sigma_2) \equiv (\sigma_1,\m_{\sigma_2}\sigma_1)$ we have:
		\begin{equation}\label{eq:HankelSign-Right}
			\operatorname{sgn}\Big(H_{n,m}\big[(\sigma_1,\m_{\sigma_2}\sigma_1)\big]\Big)  = 
			\begin{cases}
				(-1)^{nm + \frac{m(m-1)}{2}}, & n \ge m, \\[0.4em]
				(-1)^{\frac{n(n+1)}{2}}, & n \le m.
			\end{cases}
		\end{equation}
		
		\item[(ii)] For the modified Nikishin system $(\sigma_1, x \m_{\sigma_2}\sigma_1)$ we have:
		\begin{equation}\label{eq:HankelSign-mod1}
			\operatorname{sgn}\Big(H_{n,m}\big[(\sigma_1,x\m_{\sigma_2}\sigma_1)\big]\Big)  = 
			\begin{cases}
				(-1)^{nm + \frac{m(m+1)}{2}}, & n \ge m, \\[0.4em]
				(-1)^{\frac{n(n-1)}{2}}, & n \le m.
			\end{cases}
		\end{equation}
		
		\item[(iii)] For the modified Nikishin system $(x\sigma_1,  \m_{\sigma_2}\sigma_1)$ we have:
		\begin{equation}\label{eq:HankelSign-mod2}
			\operatorname{sgn}\Big(H_{n,m}\big[(x\sigma_1,\m_{\sigma_2}\sigma_1)\big]\Big)  = 
			\begin{cases}
				1, & m = 0, \\[0.4em]
				(-1)^{1 + nm + \frac{m(m+1)}{2}}, & n \ge m-2 \text{ and } m \ge 1, \\[0.4em]
				(-1)^{\frac{n(n-1)}{2}}, & n \le m-2.
			\end{cases}
		\end{equation}
	\end{itemize}
\end{proposition}
\begin{proof}
	One can obtain~\cite{KNik,Kui} the statement in \textit{(i)} by using~\eqref{eq:HankelIntegralNik} together with the generalized Andr\'eief identity~\cite{KVMLOPUC1}, the Cauchy--Vandermonde determinant~\eqref{eq:C}, and Proposition~\ref{prop:nnrr-signs}\textit{(i)}. Part \textit{(ii)} follows from Theorem~\ref{thm:XM1_Properties}\textit{(i)} combined with \textit{(i)}, and \textit{(iii)} follows from the proof of Theorem~\ref{thm:2-xm_Intro} together with \textit{(ii)}.
\end{proof}

The signs in Proposition~\ref{prop:Hsigns} immediately determine the signatures of the nearest-neighbor recurrence coefficients. Note how applying the partial Christoffel transform to the second or first measure shifts the critical line separating the two sign regimes from $n=m$ to $n=m+1$ or $n=m-1$, respectively.

\begin{theorem}
	Suppose that $\supp\, \sigma_1 \subseteq \bbR_+$ and $\supp\, \sigma_2 \subseteq \bbR_-$. The following identities hold for any $(n,m) \in \bbN^2$ and $j \in \{1, 2\}$ (with $n \ge 1$ for $a_{n,m;1}$ and $m \ge 1$ for $a_{n,m;2}$):
	\begin{itemize}
		\item[(i)] For the Nikishin system $\mathcal{N}(\sigma_1,\sigma_2) \equiv (\sigma_1,\m_{\sigma_2}\sigma_1)$ we have 
		\begin{alignat*}{3}
			\operatorname{sgn}(a_{n,m;j}) &= (-1)^{j+1}, &\quad \operatorname{sgn}(b_{n,m;2} - b_{n,m;1}) &= -1, &\quad &\text{if } n \ge m, \\
			\operatorname{sgn}(a_{n,m;j}) &= (-1)^{j},   &\quad \operatorname{sgn}(b_{n,m;2} - b_{n,m;1}) &= 1,  &\quad &\text{if } n \le m-1.
		\end{alignat*}
		
		\item[(ii)] For the modified Nikishin system $(\sigma_1, x \m_{\sigma_2}\sigma_1)$ we have:
		\begin{alignat*}{3}
			\operatorname{sgn}(a_{n,m;j}) &= (-1)^{j+1}, &\quad \operatorname{sgn}(b_{n,m;2} - b_{n,m;1}) &= -1, &\quad &\text{if } n \ge m+1, \\
			\operatorname{sgn}(a_{n,m;j}) &= (-1)^{j},   &\quad \operatorname{sgn}(b_{n,m;2} - b_{n,m;1}) &= 1,  &\quad &\text{if } n \le m.
		\end{alignat*}
		
		\item[(iii)] For the modified Nikishin system $(x\sigma_1,  \m_{\sigma_2}\sigma_1)$ we have:
		\begin{alignat*}{3}
			\operatorname{sgn}(a_{n,m;j}) &= (-1)^{j+1}, &\quad \operatorname{sgn}(b_{n,m;2} - b_{n,m;1}) &= -1, &\quad &\text{if } n \ge m-1, \\
			\operatorname{sgn}(a_{n,m;j}) &= (-1)^{j},   &\quad \operatorname{sgn}(b_{n,m;2} - b_{n,m;1}) &= 1,  &\quad &\text{if } n \le m-2,
		\end{alignat*}
		except for $m=1, j=2$, where $\operatorname{sgn}(a_{n,1;2}) = 1$ for all $n \ge 0$.
	\end{itemize}
\end{theorem}

\begin{proof}
	These results follow directly from Proposition~\ref{prop:Hsigns} and identities~\eqref{eq:a-coeffs} and~\eqref{eq:b-diff}. 
	For the standard Nikishin system \textit{(i)}, this was also established in~\cite[Thm.~1.4.1]{DenYat} using a different method (note that~\cite{DenYat} states the result for $I_1<I_2$). 
\end{proof}

\section{Examples}\label{ss:Examples}

Note that Theorems~\ref{thm:1-xm_Intro} and~\ref{thm:2-xm_Intro} (as well as the sign formulas in Section~\ref{ss:Sign}) extend to positive measures $\sigma_2$ on $\bbR_-$ of infinite total mass (such that $\int_{\bbR_-} (1+|t|)^{-1} d\sigma_2(t) < \infty$), provided $\mathfrak{m}_{\sigma_2}\sigma_1$ has finite moments: in Lemma~\ref{lem:transforms}\textit{(iii)}--\textit{(iv)} one simply has $\mu_0^{-1} = 0$, while in~\eqref{eq:jpdf} one can truncate $\sigma_2$ to $[-\Lambda, -\epsilon]$ and pass to the limit $\Lambda \to \infty$, $\epsilon \to 0^+$ by monotone convergence (cf.~\cite{FidLop11b}).

\begin{example}[Multiple Laguerre of the first kind]
	Consider $\bm{\mu} = (\mu_1, \mu_2)$ on $I_1 = [0, \infty)$ with $d\mu_j(x) = x^{\alpha_j} e^{-x} dx$, where $\alpha_1, \alpha_2 > -1$ and $\alpha_1 - \alpha_2 \notin \bbZ$. For $0 < \gamma < 1$ and $x > 0$, Euler's beta integral gives
	\begin{equation*}
		x^{-\gamma} = \mathfrak{m}_{\sigma_2}(x), \qquad d\sigma_2(t) = \frac{\sin(\pi \gamma)}{\pi} (-t)^{-\gamma} \bdone_{(-\infty,0)}(t) dt.
	\end{equation*}
	
	When $0 < \alpha_1 - \alpha_2 < 1$, setting $\gamma = \alpha_1 - \alpha_2$ and $\sigma_1 = \mu_1$ realizes $(\mu_1, \mu_2)$ as the Nikishin system $(\sigma_1, \mathfrak{m}_{\sigma_2}\sigma_1)$.
	
	When $-1 < \alpha_1 - \alpha_2 < 0$, setting $\gamma = \alpha_1 - \alpha_2 + 1 \in (0,1)$ and $\sigma_1 = \mu_1$ realizes $(\mu_1, \mu_2)$ as the system $(\sigma_1, x \mathfrak{m}_{\sigma_2}\sigma_1)$ of Theorem~\ref{thm:1-xm_Intro} at $x_0 = 0$.
	
	When $1 < \alpha_1 - \alpha_2 < 2$, setting $\gamma = \alpha_1 - \alpha_2 - 1 \in (0,1)$ and $d\sigma_1(x) = x^{\alpha_1 - 1} e^{-x} dx$ realizes $(\mu_1, \mu_2)$ as the system $(x\sigma_1, \mathfrak{m}_{\sigma_2}\sigma_1)$ of Theorem~\ref{thm:2-xm_Intro} at $x_0 = 0$. 
	
	The Jacobi--Pi\~neiro system on $[0,1]$ with weights $x^{\alpha_j}(1-x)^\beta$ can be treated in the same way.
\end{example}

\begin{example}[Multiple Meixner of the second kind]
	Consider the discrete measures on $\bbN$ with multiple Meixner weights of the second kind~\cite{Ismail},
	\begin{equation*}
		w^{(\beta, q)}(x) = \frac{(\beta)_x q^x}{x!}, \qquad \beta > 0, \quad 0 < q < 1.
	\end{equation*}
	When $0 < \beta_1 - \beta_2 < 1$, the partial fraction expansion of $\Gamma(x+\beta_2)/\Gamma(x+\beta_1)$ gives $\frac{w^{(\beta_2,q)}(x)}{w^{(\beta_1,q)}(x)} = \mathfrak{m}_{\sigma_2}(x)$, where $\sigma_2$ is a positive discrete measure supported on $\{-\beta_2 - k : k \in \bbN\} \subset (-\infty, -\beta_2]$. Here $I_1 = [0, \infty)$ and $I_2 = (-\infty, -\beta_2]$ are separated by the closed interval $\Delta = [-\beta_2, 0]$. 
	
	By Theorems~\ref{thm:1-xm_Intro} and~\ref{thm:2-xm_Intro}, for any $x_0 \in [-\beta_2, 0]$ the partially modified systems
	\begin{equation*}
		\big( w^{(\beta_1, q)}(x),\, (x - x_0)w^{(\beta_2, q)}(x) \big) \quad \text{and} \quad \big( (x - x_0)w^{(\beta_1, q)}(x),\, w^{(\beta_2, q)}(x) \big)
	\end{equation*}
	are perfect. The two endpoints of $\Delta = [-\beta_2, 0]$ yield notable special cases:
	\begin{itemize}
		\item[(a)] At $x_0 = -\beta_2$, the identity $(x + \beta_2)w^{(\beta_2, q)}(x) = \beta_2 w^{(\beta_2+1, q)}(x)$ shows that the first system reduces (up to a constant factor) to $\big(w^{(\beta_1, q)}(x), w^{(\beta_2+1, q)}(x)\big)$, which produces a parameter shift $\beta_2 \mapsto \beta_2 + 1$ like in the previous example.
		\item[(b)] At $x_0 = 0$, the identity $x w^{(\beta, q)}(x) = \beta q \, w^{(\beta+1, q)}(x-1) \bdone_{\{x \ge 1\}}$ shows that
		\begin{equation*}
			\big( w^{(\beta_1, q)}(x),\, w^{(\beta_2+1, q)}(x-1)\bdone_{\{x \ge 1\}} \big) \quad \text{and} \quad \big( w^{(\beta_1+1, q)}(x-1)\bdone_{\{x \ge 1\}},\, w^{(\beta_2, q)}(x) \big)
		\end{equation*}
		are perfect systems on $\bbN$. Note that these are the original Meixner weights with one of the parameters shifted by $1$ and its support shifted to $\bbN + 1$.
	\end{itemize}
	Multiple Hahn weights on $\{0, 1, \dots, M\}$ can be treated in the same way.
\end{example}


\bibsection

\end{document}